\documentclass[12pt]{article}

\usepackage{amsmath, epsfig, cite}
\usepackage{amssymb}
\usepackage{amsfonts}
\usepackage{latexsym}
\usepackage{float}
\usepackage{mathrsfs}
\usepackage{booktabs}
\usepackage{caption}

\newtheorem{thm}{Theorem}[section]

\newtheorem{cor}[thm]{Corollary}
\newtheorem{conj}[thm]{Conjecture}
\newtheorem{lem}[thm]{Lemma}

\numberwithin{equation}{section}

\newcommand{\qed}{{\hfill$\square$}\medskip}
\UseRawInputEncoding
\begin{document}


\begin{center}
{\large\bf A general positivity result on coefficients of certain $q$-series}
\end{center}

\vskip 2mm \centerline{Ji-Cai Liu}
\begin{center}
{\footnotesize Department of Mathematics, Wenzhou University, Wenzhou 325035, PR China\\
{\tt jcliu2016@gmail.com} }
\end{center}


\vskip 0.7cm \noindent{\bf Abstract.}
Based on a classical result on partitions of an integer into a finite set of positive integers,
we establish a general positivity result on coefficients of certain $q$-series which uniformly refines the positivity of truncated pentagonal number series, truncated Gauss' identities and some special cases of truncated Jacobi triple product identity. As an application, we prove two positivity conjectures due to Merca.


\vskip 2mm
\noindent{\it MR Subject Classifications}: 05A17, 05A20
\section{Introduction}
One of Euler's most profound discoveries is the pentagonal number theorem (see \cite[Corollary 1.7]{andrews-b-1998}):
\begin{align}
(q;q)_{\infty}=\sum_{j=-\infty}^{\infty}(-1)^jq^{j(3j+1)/2}.\label{a-euler}
\end{align}
Here and throughout this paper, the $q$-shifted factorial is defined by $(a;q)_0=1$,
$(a;q)_n=(1-a)(1-aq)\cdots (1-aq^{n-1})$ for $n\ge 1$, $(a;q)_{\infty}=\prod_{k=0}^{\infty}(1-aq^k)$ and
$(a_1,a_2,\cdots,a_m;q)_{\infty}=(a_1;q)_{\infty}(a_2;q)_{\infty}\cdots (a_m;q)_{\infty}$.
The $q$-binomial coefficient is defined as
\begin{align*}
{n\brack k}={n\brack k}_{q}
=\begin{cases}
\displaystyle\frac{(q;q)_n}{(q;q)_k(q;q)_{n-k}} &\text{if $0\le k\le n$},\\[10pt]
0 &\text{otherwise.}
\end{cases}
\end{align*}

Let $p(n)$ denote the number of partitions of $n$. The partition function $p(n)$ has the generating function:
\begin{align*}
\sum_{n=0}^{\infty}p(n)q^n=\frac{1}{(q;q)_{\infty}}.
\end{align*}

Andrews and Merca \cite{am-jcta-2012} showed that for $k\ge 1$,
\begin{align*}
(-1)^{k-1}\sum_{j=0}^{k-1}(-1)^j\left(p(n-j(3j+1)/2)-p(n-(j+1)(3j+2)/2)\right)\ge 0,
\end{align*}
which is equivalent to
\begin{align}
\frac{(-1)^{k-1}}{(q;q)_{\infty}}\sum_{j=0}^{k-1} (-1)^jq^{j(3j+1)/2}(1-q^{2j+1})+(-1)^k
\in \mathbb{N}[[q]].\label{a-2}
\end{align}
In order to prove \eqref{a-2}, Andrews and Merca \cite{am-jcta-2012} established the truncation of the pentagonal number theorem \eqref{a-euler}:
\begin{align*}
\frac{1}{(q;q)_{\infty}}\sum_{j=0}^{k-1} (-1)^jq^{j(3j+1)/2}(1-q^{2j+1})
=1+(-1)^{k-1}\sum_{j=k}^{\infty}\frac{q^{k(k-1)/2+(k+1)j}}{(q;q)_j}{j-1\brack k-1}.
\end{align*}

By \eqref{a-euler}, we have
\begin{align}
\sum_{j=0}^{k-1} (-1)^jq^{j(3j+1)/2}(1-q^{2j+1})
=(q;q)_{\infty}-\sum_{j\not \in [-k,k-1]}(-1)^jq^{j(3j+1)/2}.\label{a-3}
\end{align}
Here and throughout this paper, we use the notation: for integers $a$ and $b$ with $a<b$,
\begin{align*}
\sum_{j\not\in[a,b]} A_j=\sum_{j=-\infty}^{a-1}A_j+\sum_{j=b+1}^{\infty}A_j.
\end{align*}
From \eqref{a-3}, we deduce that \eqref{a-2} is equivalent to
\begin{align}
\frac{1}{(q;q)_{\infty}}\sum_{j\not \in [-k,k-1]}(-1)^{j+k}q^{j(3j+1)/2}
\in \mathbb{N}[[q]].\label{a-4}
\end{align}

Motivated by the work of Andrews and Merca \cite{am-jcta-2012}, Guo and Zeng \cite{gz-jcta-2013} investigated two truncated identities of Gauss (see \cite[page 23]{andrews-b-1998}):
\begin{align}
\sum_{j=-\infty}^{\infty}(-1)^jq^{j^2}=\frac{(q;q)_{\infty}}{(-q;q)_{\infty}},\label{a-5}
\end{align}
and
\begin{align}
\sum_{j=0}^{\infty}(-1)^jq^{j(2j+1)}(1-q^{2j+1})
=\frac{(q^2;q^2)_{\infty}}{(-q;q^2)_{\infty}}.\label{a-6}
\end{align}
Guo and Zeng \cite{gz-jcta-2013} showed that for $k\ge 1$,
\begin{align}
\frac{(-q;q)_{\infty}}{(q;q)_{\infty}}\sum_{j=1-k}^{k-1} (-1)^jq^{j^2}=1+(-1)^{k-1}\sum_{j=k}^{\infty}
\frac{(-q;q)_{k-1}(-1;q)_{j-k+1}q^{jk}}{(q;q)_j}{j-1\brack k-1},\label{a-7}
\end{align}
and
\begin{align}
&\frac{(-q;q^2)_{\infty}}{(q^2;q^2)_{\infty}}\sum_{j=0}^{k-1}(-1)^jq^{j(2j+1)}(1-q^{2j+1})\notag\\[7pt]
&=1+(-1)^{k-1}\sum_{j=k}^{\infty}\frac{(-q;q^2)_k (-q;q^2)_{j-k}q^{2(k+1)j-k}}{(q^2;q^2)_j}{j-1\brack k-1}_{q^2}.\label{a-8}
\end{align}

Note that
\begin{align*}
\sum_{n=0}^{\infty}\overline{p}(n)q^n&=\frac{(-q;q)_{\infty}}{(q;q)_{\infty}},\\[7pt]
\sum_{n=0}^{\infty}\text{pod}(n)q^n&=\frac{(-q;q^2)_{\infty}}{(q^2;q^2)_{\infty}},
\end{align*}
where the overpartition function $\overline{p}(n)$ denotes the number of ways of writing the integer $n$ as a sum of positive integers in non-increasing order in which the first occurrence of an integer may be
overlined or not (see \cite{cl-tams-2004}), and $\text{pod}(n)$ denotes the number of partitions of $n$ wherein odd parts are distinct (see \cite{hs-rj-2010}).

From \eqref{a-7} and \eqref{a-8}, Guo and Zeng \cite{gz-jcta-2013} deduced that for $n,k\ge 1$,
\begin{align}
(-1)^k\left(\overline{p}(n)+2\sum_{j=1}^k(-1)^j\overline{p}(n-j^2)\right)\ge 0,\label{a-9}
\end{align}
and
\begin{align}
(-1)^{k-1}\sum_{j=0}^{k-1}(-1)^j\left(\text{pod}(n-j(2j+1))-\text{pod}(n-(j+1)(2j+1))\right)\ge 0.\label{a-10}
\end{align}
The result \eqref{a-9} was strengthened by Mao \cite{mao-jcta-2015}, Yee \cite{yee-jcta-2015} and Wang--Yee \cite{wy-am-2020} in different approaches as follows:
\begin{align}
(-1)^{k-1}\left(\overline{p}(n)+2\sum_{j=1}^{k-1}(-1)^j\overline{p}(n-j^2)\right)-\overline{p}(n-k^2)\ge 0,\label{a-11}
\end{align}
which was originally conjectured by Guo and Zeng \cite{gz-jcta-2013}.

Through the same discussion as above, we find that \eqref{a-10} and \eqref{a-11} are equivalent to
\begin{align}
\frac{(-q;q)_{\infty}}{(q;q)_{\infty}}\sum_{j\not\in [-k,k-1]}(-1)^{j+k} q^{j^2}
\in \mathbb{N}[[q]],\label{a-12}
\end{align}
and
\begin{align}
\frac{(-q;q^2)_{\infty}}{(q^2;q^2)_{\infty}}\sum_{j\not\in[-k,k-1]}(-1)^{j+k}q^{2j^2+j}\in \mathbb{N}[[q]].\label{a-13}
\end{align}

The Jacobi triple product identity \cite[page 49]{andrews-b-1998} implies that
\begin{align*}
(q,q^4,q^5;q^5)_{\infty}=\sum_{j=-\infty}^{\infty}(-1)^jq^{j(5j+3)/2},
\end{align*}
and
\begin{align*}
(q^2,q^3,q^5;q^5)_{\infty}=\sum_{j=-\infty}^{\infty}(-1)^jq^{j(5j+1)/2}.
\end{align*}
Guo and Zeng \cite{gz-jcta-2013} conjectured that for $k\ge1$,
\begin{align}
(-1)^k+\frac{(-1)^{k-1}}{(q,q^4,q^5;q^5)_{\infty}}\sum_{j=-k}^{k-1}(-1)^jq^{j(5j+3)/2}\in \mathbb{N}[[q]],\label{a-14}
\end{align}
and
\begin{align}
(-1)^k+\frac{(-1)^{k-1}}{(q^2,q^3,q^5;q^5)_{\infty}}\sum_{j=-k}^{k-1}(-1)^jq^{j(5j+1)/2}\in \mathbb{N}[[q]],\label{a-15}
\end{align}
which were proved by Mao \cite{mao-jcta-2015} and Yee \cite{yee-jcta-2015} in different methods.
Through the same discussion as above, we find that \eqref{a-14} and \eqref{a-15} are equivalent to
\begin{align}
\frac{1}{(q,q^4,q^5;q^5)_{\infty}}\sum_{j\not\in[-k,k-1]}(-1)^{j+k}q^{j(5j+3)/2}
\in \mathbb{N}[[q]],\label{a-16}
\end{align}
and
\begin{align}
\frac{1}{(q^2,q^3,q^5;q^5)_{\infty}}\sum_{j\not\in[-k,k-1]}(-1)^{j+k}q^{j(5j+1)/2}
\in \mathbb{N}[[q]].\label{a-17}
\end{align}

Recently, Yao \cite{yao-jcta-2024} strengthened \eqref{a-4} as follows: for $k\ge 1$ and $(a,b,c)=(1,2,3)$,
\begin{align}
\frac{1}{(1-q^a)(1-q^b)(1-q^c)}\sum_{j\not\in [-k,k-1]} (-1)^{j+k}q^{j(3j+1)/2}
\in \mathbb{N}[[q]].\label{a-18}
\end{align}
Subsequently, Zhou \cite{zhou-jcta-2024} found the following more triples $(a,b,c)$ such that \eqref{a-18} holds for $k\ge 1$:
\begin{align*}
(1,2,3),(1,2,5),(1,2,7),(1,3,4),(1,3,5).
\end{align*}
Note that
\begin{align*}
\frac{1}{(q;q)_{\infty}}=\frac{P(q)}{(1-q^a)(1-q^b)(1-q^c)},
\end{align*}
where $P(q)\in \mathbb{N}[[q]]$ and $(a,b,c)\in\left\{(1,2,3),(1,2,5),(1,2,7),(1,3,4),(1,3,5)\right\}$.
The results due to Yao \cite{yao-jcta-2024} and Zhou \cite{zhou-jcta-2024} are stronger than \eqref{a-4}.

It is natural to consider Yao-Zhou type extensions of \eqref{a-12}, \eqref{a-13}, \eqref{a-16} and \eqref{a-17}. The motivation of the paper is to establish a general positivity result which extends
\eqref{a-4}, \eqref{a-12}, \eqref{a-13}, \eqref{a-16} and \eqref{a-17} uniformly.

Let $Ax^2+Bx$ be an integer valued polynomial with $A> B\ge 0$, and $a,b,c$ be distinct positive integers with $(a,b)=(a,c)=(b,c)=1$.
For integers $k\ge 1$ and $n\ge 0$, the coefficient $\gamma_{a,b,c,A,B}^k (n)$ is defined by
\begin{align}
\frac{1}{(1-q^a)(1-q^b)(1-q^c)}\sum_{j\not\in [-k,k-1]} (-1)^{j+k}q^{Aj^2+Bj}=\sum_{n= 0}^{\infty}\gamma_{a,b,c,A,B}^k (n)q^n.\label{a-gamma}
\end{align}

What we want to do is to determine the values $K_{a,b,c,A,B}$ and $N_{a,b,c,A,B}^k$ such that
\begin{align*}
\gamma_{a,b,c,A,B}^k(n)\ge 0\quad\text{for $n\ge 0$ with $k\ge K_{a,b,c,A,B}$},
\end{align*}
and
\begin{align*}
\gamma_{a,b,c,A,B}^k(n)\ge 0\quad\text{for $n\ge N_{a,b,c,A,B}^k$}.
\end{align*}
Once the values $K_{a,b,c,A,B}$ and $N_{a,b,c,A,B}^k$ are determined, the remaining task is to verify a finite number of values $\gamma_{a,b,c,A,B}^k (n)$ for $0\le n<N_{a,b,c,A,B}^k$ with $1\le k<K_{a,b,c,A,B}$ through mathematical software such as Maple. The proof of the main theorem is inspired by Zhou's method \cite{zhou-jcta-2024}, in which a classical result on partitions of an integer into a finite set of positive integers plays an important role (see Lemma \ref{b-lemma}).

As an application, we prove two positivity conjectures due to Merca \cite[Conjectures 13 and 15]{merca-am-2024}:
\begin{conj}[Merca]
For a positive integer $n$, let $\nu_2(n)$ denote the $2$-adic order of $n$ and $N_n=n(1+\nu_2(n)/2)$.
For all positive integers $k$, we have
\begin{align}
(-1)^k\left(1-\frac{1}{(q;q)_{\infty}}\sum_{j=1-k}^k (-1)^jq^{j(3j-1)/2}\right)\prod_{n=1}^{\infty}(1-q^{N_{2n}})\in \mathbb{N}[[q]],\label{merca-conj-1}
\end{align}
and
\begin{align}
(-1)^{k-1}\left(1-\frac{1}{(q;q)_{\infty}}\sum_{j=-k}^k (-1)^jq^{j(3j-1)/2}\right)\prod_{n=1}^{\infty}(1-q^{N_{2n}})\in \mathbb{N}[[q]].\label{merca-conj-2}
\end{align}
\end{conj}

The rest of the paper is organized as follows. The main results are stated in the next section. Section 3 is devoted to the proof of Lemma \ref{b-lemma}. The proof of the main theorem is presented in Section 4. We prove \eqref{merca-conj-1} and \eqref{merca-conj-2} in the last section.

\section{Main results}
In order to state the main theorem, we require some notation.
Let $G(x)=ux^2+vx+w$ be a real quadratic function in variable $x$ with $u>0$ and $\Delta_G=v^2-4uw$. The operation $\mathcal{R}_x$ is defined as
\begin{align*}
\mathcal{R}_xG(x)=
\begin{cases}
1&\quad\text{if $\Delta_G<0$,}\\[10pt]
\displaystyle \max\left\{1,\frac{-v+\sqrt{\Delta_G}}{2u}\right\}&\quad\text{if $\Delta_G\ge 0$.}
\end{cases}
\end{align*}
It is clear that $G(x)\ge 0$ for all real numbers $x\ge \mathcal{R}_xG(x)$.
For real quadratic functions $G_i(x)=u_ix^2+v_ix+w_i$ with $u_i>0$ for $i=1,2,\cdots,n$, let
\begin{align*}
\mathcal{T}_{x}\left\{G_1(x),G_2(x),\cdots,G_n(x)\right\}
=\max\left\{\mathcal{R}_xG_1(x),\mathcal{R}_xG_2(x),\cdots,\mathcal{R}_xG_n(x)\right\}.
\end{align*}
It is also clear that $G_i(x)\ge 0$ for $i=1,\cdots,n$ and all real numbers
\begin{align*}
x\ge \mathcal{T}_{x}\left\{G_1(x),G_2(x),\cdots,G_n(x)\right\}.
\end{align*}

We also need the notation $D_{a,b,c}$ in the main theorem, which is related to the following result due to P\'olya and Szeg\H{o} \cite[Problem 27.1, page 5]{ps-b-1998}.
\begin{lem}\label{b-lemma}
For distinct positive integers $a,b,c$ with $(a,b)=(a,c)=(b,c)=1$, let
\begin{align*}
\frac{1}{(1-q^a)(1-q^b)(1-q^c)}=\sum_{n=0}^{\infty}\alpha(n)q^n,
\end{align*}
and
\begin{align*}
\beta(n)=\alpha(n)-\frac{n^2+(a+b+c)n}{2abc}.
\end{align*}
Then $\{\beta(n)\}_{n\ge 0}$ has a period of $abc$.
\end{lem}

We remark that P\'olya and Szeg\H{o} \cite[Problem 27.1, page 5]{ps-b-1998} qualitatively described $F(x)=(x^2+(a+b+c)x)/(2abc)$ as a polynomial with rational coefficients of degree $2$. Zhou \cite{zhou-jcta-2024} gave the explicit expression of $F(x)$ without proof. We shall present a complete proof of Lemma \ref{b-lemma} in Section 3.

Since $\{\beta(n)\}_{n\ge 0}$ is a periodic sequence, there exists a smallest bound $D_{a,b,c}$ such that
$|\beta(n)|\le D_{a,b,c}$ for all integers $n\ge 0$.

For a real number $x$, let $\lceil x\rceil$ denote the smallest integer greater than or equal to $x$.
Now we are ready to state the main theorem.

\begin{thm}\label{t-1}
Let $\gamma_{a,b,c,A,B}^k (n)$ be defined by \eqref{a-gamma}.
\begin{itemize}
\item[(1)] For all $n\ge 0$ with $k\ge K_{a,b,c,A,B}$, we have $\gamma_{a,b,c,A,B}^k (n)\ge 0$, where $K_{a,b,c,A,B}$ is given by
    \begin{align*}
    K_{a,b,c,A,B}=\lceil\mathcal{T}_{k}\left\{H_1,H_2,\cdots,H_{13}\right\}\rceil
    \end{align*}
    with
{\scriptsize\begin{align*}
H_1&=\frac{2(A-B)^2k^2}{abc}+\frac{(A-B)(2A-2B+a+b+c)k}{abc}+\frac{(A-B)(A-B+a+b+c)}{2abc}-2D_{a,b,c},\\[5pt]
H_2&=\frac{2(A-B)(A+B)k^2}{abc}+\frac{(A-B)(2A+4B+a+b+c)k}{abc}+\frac{(A-B)(A+3B+a+b+c)}{2abc}-3D_{a,b,c},\\[5pt]
H_3&=\frac{4A(A-B)k^2}{abc}+\frac{2(2A+B)(A-B)k}{abc}-\frac{(A-B)(A-3B+a+b+c)}{abc}-3D_{a,b,c},\\[5pt]
H_4&=\frac{4A(A-B)k^2}{abc}+\frac{2(A-B)(2A-B)k}{abc}-\frac{(A-B)(A+5B+a+b+c)}{abc}-4D_{a,b,c},\\[5pt]
H_5&=\frac{4A(A-B)k^2}{abc}+\frac{2(A-B)(2A-B)k}{abc}-\frac{24(A-B)(A+B)+(a+b+c)^2}{8abc}-5D_{a,b,c},\\[5pt]
H_6&=\frac{4A(A-B)k^2}{abc}+\frac{6(A-B)(2A-B)k}{abc}+\frac{(A-B)(A-11B+a+b+c)}{abc}-4D_{a,b,c},\\[5pt]
H_7&=\frac{2(A-B)(3A-B)k^2}{abc}+\frac{(A-B)(10A-8B+a+b+c)k}{abc}+\frac{3(A-B)(A-5B+a+b+c)}{2abc}-6D_{a,b,c},\\[5pt]
H_8&=\frac{4A(A-B)k^2}{abc}+\frac{6(A-B)(2A+B)k}{abc}+\frac{(A-B)(A+13B+a+b+c)}{abc}-4D_{a,b,c},\\[5pt]
H_9&=\frac{2(A-B)(3A+B)k^2}{abc}+\frac{(A-B)(10A+a+b+c+10B)k}{abc}+\frac{3(A-B)(A+7B+a+b+c)}{2abc}-7D_{a,b,c},\\[5pt]
H_{10}&=\frac{4A(A-B)k^2}{abc}+\frac{2(A-B)(8A+B)k}{abc}+\frac{(A-B)(A-15B+a+b+c)}{abc}
-4D_{a,b,c},\\[5pt]
H_{11}&=\frac{8A(A-B)k^2}{abc}+\frac{4(A-B)(4A+B)k}{abc}-\frac{2(A-B)(A-7B+a+b+c)}{abc}-7D_{a,b,c},\\[5pt]
H_{12}&=\frac{4A(A-B)k^2}{abc}+\frac{2(A-B)(8A-B)k}{abc}-\frac{(A-B)(A+17B+a+b+c)}{abc}-4D_{a,b,c},\\[5pt]
H_{13}&=\frac{8A(A-B)k^2}{abc}+\frac{4(A-B)(4A-B)k}{abc}-\frac{2(A-B)(A+9B+a+b+c)}{abc}-8D_{a,b,c}.
\end{align*}}

\item[(2)] For all $n\ge N_{a,b,c,A,B}^k=A\left(k+2L_{a,b,c,A,B}^k\right)^2+B\left(k+2L_{a,b,c,A,B}^k\right)$, we have $\gamma_{a,b,c,A,B}^k (n)\ge 0$, where $L_{a,b,c,A,B}^k$ is given by
    \begin{align*}
    L_{a,b,c,A,B}^k=\lceil\mathcal{T}_{l}\left\{G_1,G_2,\cdots,G_5\right\}\rceil
    \end{align*}
    with
    {\scriptsize\begin{align*}
    G_1&=\frac{2(A-B)(2Ak-A-B)l^2}{abc}+\frac{(2k+1)(A-B)(2Ak-A-B)l}{abc}-\frac{(a+b+c)^2}{8abc}-(4l+1)D_{a,b,c},\\[5pt]
    G_2&=\frac{4(A-B)(Ak-B)l^2}{abc}+\frac{(A-B)(4Ak^2+4Ak-6Bk+A-3B+a+b+c)l}{abc}\\[5pt]
    &+\frac{(2k+1)(A-B)(2Ak-2Bk+A-B+a+b+c)}{2abc}-(4l+2)D_{a,b,c},\\[5pt]
    G_3&=\frac{4(A-B)(Ak+B)l^2}{abc}+\frac{(A-B)(4Ak^2+4Ak+6Bk+A+5B+a+b+c)l}{abc}\\[5pt]
    &+\frac{(2k+1)(A-B)(2Ak+2Bk+A+3B+a+b+c)}{2abc}-(4l+3)D_{a,b,c},\\[5pt]
    G_4&=\frac{4(A-B)(Ak+B)l^2}{abc}+\frac{(A-B)(4Ak^2+8Ak+2Bk-A+7B-a-b-c)l}{abc}\\[5pt]
    &+\frac{(A-B)(4Ak^2+4Ak+2Bk-A+3B-a-b-c)}{abc}-(4l+3)D_{a,b,c},\\[5pt]
    G_5&=\frac{4(A-B)(Ak-B)l^2}{abc}+\frac{(A-B)(4Ak^2+8Ak-2Bk-A-9B-a-b-c)l}{abc}\\[5pt]
       &+\frac{(A-B)(4Ak^2+4Ak-2Bk-A-5B-a-b-c)}{abc}-(4l+4)D_{a,b,c}.
    \end{align*}}
\end{itemize}
\end{thm}

From Theorem \ref{t-1}, we derive several corollaries which extend \eqref{a-4}, \eqref{a-12}, \eqref{a-13}, \eqref{a-16} and \eqref{a-17}.
\begin{cor}\label{cor-3}
For all integers $k\ge 1,n\ge 0$ and
\begin{align*}
(a,b,c)\in\left\{(1,2,3),(1,2,5),(1,2,7),(1,3,4),(1,3,5),(1,3,8),(1,4,5),(1,4,7)\right\},
\end{align*}
we have
$\gamma_{a,b,c,3/2,1/2}^k (n)\ge 0$.
\end{cor}
{\noindent\bf Remark.}The triples $(1,3,8),(1,4,5)$ and $(1,4,7)$ were not listed by Zhou \cite{zhou-jcta-2024}.

To prove the above result, we only need calculate the values $K_{a,b,c,A,B}$ and $N_{a,b,c,A,B}^k$ with $1\le k<K_{a,b,c,A,B}$, and verify a finite number of values $\gamma_{a,b,c,A,B}^k (n)$ for $0\le n<N_{a,b,c,A,B}^k$ with $1\le k<K_{a,b,c,A,B}$.

\begin{table}[H]
\caption*{\text{The case $A=3/2$ and $B=1/2$}}
\centering
\begin{tabular}{cccc}
\toprule
$(a,b,c)$&$D_{a,b,c}$&$K_{a,b,c,A,B}$&$\left\{N_{a,b,c,A,B}^k\right\}_{k=1}^{K_{a,b,c,A,B}-1}$\\
\midrule
$(1,2,3)$&$1$&$3$&$805,57$\\[10pt]
$(1,2,5)$&$1$&$3$&$2301,155$\\[10pt]
$(1,2,7)$&$8/7$&$4$&$5985,392,126$\\[10pt]
$(1,3,4)$&$1$&$4$&$3337,222,77$\\[10pt]
$(1,3,5)$&$1$&$4$&$5251,301,77$\\[10pt]
$(1,3,8)$&$17/16$&$5$&$15352,876,260,100$\\[10pt]
$(1,4,5)$&$9/8$&$5$&$11926,737,187,100$\\[10pt]
$(1,4,7)$&$8/7$&$6$&$24257, 1365, 442, 155, 126$\\
\bottomrule
\end{tabular}
\end{table}

\begin{cor}
For all integers $k\ge 1,n\ge 0$ and
\begin{align*}
(a,b,c)\in\left\{(1,2,3),(1,2,5),(1,3,5)\right\},
\end{align*}
we have $\gamma_{a,b,c,1,0}^k (n)\ge 0$.
\end{cor}

\begin{table}[H]
\caption*{\text{The case $A=1$ and $B=0$}}
\centering
\begin{tabular}{cccc}
\toprule
$(a,b,c)$&$D_{a,b,c}$&$K_{a,b,c,A,B}$&$\left\{N_{a,b,c,A,B}^k\right\}_{k=1}^{K_{a,b,c,A,B}-1}$\\
\midrule
$(1,2,3)$&$1$&$3$&$529,64$\\[10pt]
$(1,2,5)$&$1$&$4$&$1521,144,49$\\[10pt]
$(1,3,5)$&$1$&$5$&$3481,400,121,64$\\
\bottomrule
\end{tabular}
\end{table}

\begin{cor}
For all integers $k\ge 1,n\ge 0$ and
\begin{align*}
(a,b,c)\in\left\{(1,3,4),(1,3,5),(1,4,5),(1,5,7)\right\},
\end{align*}
we have $\gamma_{a,b,c,2,1}^k (n)\ge 0$.
\end{cor}

\begin{table}[H]
\caption*{\text{The case $A=2$ and $B=1$}}
\centering
\begin{tabular}{cccc}
\toprule
$(a,b,c)$&$D_{a,b,c}$&$K_{a,b,c,A,B}$&$\left\{N_{a,b,c,A,B}^k\right\}_{k=1}^{K_{a,b,c,A,B}-1}$\\
\midrule
$(1,3,4)$&$1$&$3$&$4465,136$\\[10pt]
$(1,3,5)$&$1$&$4$&$7021,300,105$\\[10pt]
$(1,4,5)$&$9/8$&$4$&$15931,528,171$\\[10pt]
$(1,5,7)$&$1$&$5$&$38781,1378,351,136$\\
\bottomrule
\end{tabular}
\end{table}

\begin{cor}
For all integers $k\ge 1$ and $n\ge 0$,
we have $\gamma_{1,4,5,5/2,3/2}^k (n)\ge 0$.
\end{cor}

\begin{table}[H]
\caption*{\text{The case $A=5/2$ and $B=3/2$}}
\centering
\begin{tabular}{cccc}
\toprule
$(a,b,c)$&$D_{a,b,c}$&$K_{a,b,c,A,B}$&$\left\{N_{a,b,c,A,B}^k\right\}_{k=1}^{K_{a,b,c,A,B}-1}$\\
\midrule
$(1,4,5)$&$9/8$&$4$&$19936, 511, 133$\\
\bottomrule
\end{tabular}
\end{table}

\begin{cor}
For all integers $k\ge 1$ and $n\ge 0$, we have $\gamma_{2,3,5,5/2,1/2}^k (n)\ge 0$.
\end{cor}

\begin{table}[H]
\caption*{\text{The case $A=5/2$ and $B=1/2$}}
\centering
\begin{tabular}{cccc}
\toprule
$(a,b,c)$&$D_{a,b,c}$&$K_{a,b,c,A,B}$&$\left\{N_{a,b,c,A,B}^k\right\}_{k=1}^{K_{a,b,c,A,B}-1}$\\
\midrule
$(2,3,5)$&$1$&$3$&$2117,164$\\
\bottomrule
\end{tabular}
\end{table}

\section{Proof of Lemma \ref{b-lemma}}
Since $(a,b)=(a,c)=(b,c)=1$, we have $(1-q)^3,1+q+\cdots+q^{a-1},1+q+\cdots+q^{b-1}$ and $1+q+\cdots+q^{c-1}$ are pairwise coprime. By the partial fraction decomposition, we obtain
\begin{align}
&\frac{1}{(1-q^a)(1-q^b)(1-q^c)}\notag\\[5pt]
&=\frac{1}{(1-q)^3(1+q+\cdots+q^{a-1})(1+q+\cdots+q^{b-1})(1+q+\cdots+q^{c-1})}\notag\\[5pt]
&=\frac{R_a(q)}{1+q+\cdots+q^{a-1}}+\frac{R_b(q)}{1+q+\cdots+q^{b-1}}+\frac{R_c(q)}{1+q+\cdots+q^{c-1}}\notag\\[5pt]
&+\frac{r_1}{1-q}+\frac{r_2}{(1-q)^2}+\frac{r_3}{(1-q)^3},\label{d-1}
\end{align}
where $r_1,r_2$ and $r_3$ are rational numbers and $R_a(q),R_b(q)$ and $R_c(q)$ are polynomials with rational coefficients of degrees less than $a-1,b-1,c-1$, respectively.

By the L'H\^{o}pital's rule, we have
\begin{align}
r_3&=\lim_{q\to 1}\frac{(1-q)^3}{(1-q^a)(1-q^b)(1-q^c)}=\frac{1}{abc},\label{d-2}
\end{align}
and
\begin{align}
r_2=\lim_{q\to 1}(1-q)^2\left(\frac{1}{(1-q^a)(1-q^b)(1-q^c)}-\frac{r_3}{(1-q)^3}\right)=\frac{a+b+c-3}{2abc}.\label{d-3}
\end{align}
It follows from \eqref{d-2} and \eqref{d-3} that
\begin{align}
&\frac{r_2}{(1-q)^2}+\frac{r_3}{(1-q)^3}\notag\\[5pt]
&=\frac{1}{abc}\sum_{n\ge 0}{-3\choose n}(-q)^n
+\frac{a+b+c-3}{2abc}\sum_{n\ge 0}{-2\choose n}(-q)^n\notag\\[5pt]
&=\frac{1}{abc}\sum_{n\ge 0} \frac{(n+1)(n+2)}{2}q^n+\frac{a+b+c-3}{2abc}\sum_{n\ge 0}(n+1)q^n\notag\\
&=\sum_{n\ge 0}\frac{n^2+(a+b+c)n}{2abc}q^n+\sum_{n\ge 0}\frac{a+b+c-1}{2abc}q^n.\label{d-4}
\end{align}

Combining \eqref{d-1} and \eqref{d-4}, we arrive at
\begin{align}
&\frac{1}{(1-q^a)(1-q^b)(1-q^c)}\notag\\[5pt]
&=\sum_{n\ge 0}\frac{n^2+(a+b+c)n}{2abc}q^n\notag\\[5pt]
&+\sum_{n\ge 0}\frac{a+b+c-1}{2abc}q^n+\frac{r_1}{1-q}+\frac{(1-q)R_a(q)}{1-q^{a}}
+\frac{(1-q)R_b(q)}{1-q^{b}}+\frac{(1-q)R_c(q)}{1-q^{c}}.\label{d-5}
\end{align}
Let
\begin{align*}
&\sum_{n\ge 0}\frac{a+b+c-1}{2abc}q^n+\frac{r_1}{1-q}=\sum_{n\ge 0}t_1(n) q^n,\\[5pt]
&\frac{(1-q)R_a(q)}{1-q^{a}}=\sum_{n\ge 0}t_a(n)q^n,\\[5pt]
&\frac{(1-q)R_b(q)}{1-q^{b}}=\sum_{n\ge 0}t_b(n)q^n,\\[5pt]
&\frac{(1-q)R_c(q)}{1-q^{c}}=\sum_{n\ge 0}t_c(n)q^n,
\end{align*}
and
\begin{align*}
\sum_{n\ge 0}\frac{a+b+c-1}{2abc}q^n+\frac{r_1}{1-q}+\frac{(1-q)R_a(q)}{1-q^{a}}
+\frac{(1-q)R_b(q)}{1-q^{b}}+\frac{(1-q)R_c(q)}{1-q^{c}}=\sum_{n\ge 0}t(n)q^n.
\end{align*}
Noting that $\{t_1(n)\}_{n\ge 0},\{t_a(n)\}_{n\ge 0},\{t_b(n)\}_{n\ge 0}$ and $\{t_c(n)\}_{n\ge 0}$
have periods of $1,a,b$ and $c$, respectively, we conclude that $\{t(n)\}_{n\ge 0}$ has a period of $abc$.

Finally, we rewrite \eqref{d-5} as
\begin{align*}
\frac{1}{(1-q^a)(1-q^b)(1-q^c)}
=\sum_{n\ge 0}\frac{n^2+(a+b+c)n}{2abc}q^n+\sum_{n\ge 0}t(n)q^n.
\end{align*}
This completes the proof of Lemma \ref{b-lemma}.

{\noindent\bf Remark.} By using the same method as in the proof of Lemma \ref{b-lemma}, we can also show that for pairwise coprime positive integers $a,b,c,d$ and $e$,
\begin{align}
&\frac{1}{(1-q^a)(1-q^b)(1-q^c)(1-q^d)}\notag\\[5pt]
&=\sum_{n\ge 0}\frac{2n^3+3(a+b+c+d)n^2+\left(a^2+b^2+c^2+d^2+3(ab+ac+ad+bc+bd+cd)\right)n}{12abcd}q^n\notag\\[5pt]
&+\sum_{n\ge 0}t_{a,b,c,d}(n)q^n,\label{d-6}
\end{align}
and
\begin{align}
&\frac{1}{(1-q^a)(1-q^b)(1-q^c)(1-q^d)(1-q^e)}\notag\\[5pt]
&=\sum_{n\ge 0}\frac{n^4+2C_1n^3+(C_2+3C_3)n^2+C_1C_3 n}{24abcde}q^n+\sum_{n\ge 0}t_{a,b,c,d,e}(n)q^n,\label{d-7}
\end{align}
where
\begin{align*}
&C_1=a+b+c+d+e,\\[5pt]
&C_2=a^2+b^2+c^2+d^2+e^2,\\[5pt]
&C_3=ab+ac+ad+ae+bc+bd+be+cd+ce+de,
\end{align*}
and $\{t_{a,b,c,d}(n)\}_{n\ge 0}$ and $\{t_{a,b,c,d,e}(n)\}_{n\ge 0}$ have periods of $abcd$ and $abcde$, respectively.

Note that \eqref{d-7} is the general form of the results due to Chen and Yao \cite[Lemmas 2.2 and 3.1]{cy-jcta-2024}.

\section{Proof of Theorem \ref{t-1}}
In order to prove Theorem \ref{t-1}, we require a trivial result.
\begin{lem}\label{c-lemma}
Let $Y(x)=ux^2+vx+w$ be a real quadratic function with $u>0$.
\begin{itemize}
\item[(1)] If $2u+v\ge 0$ and $u+v+w\ge 0$, then
$Y(x)$ has real roots $x_1$ and $x_2$ with $x_1,x_2\le 1$ or $Y(x)$ has no real root.
\item[(2)] If $w\le 0$ and $u+v+w\ge 0$, then
$Y(x)$ has real roots $x_1$ and $x_2$ with $x_1,x_2\le 1$ or $Y(x)$ has no real root.
\end{itemize}
\end{lem}

Now we are ready to prove Theorem \ref{t-1}.

Let $f(j)=Aj^2+Bj,~~g(j)=Aj^2-Bj$ and
\begin{align*}
F(x)=\frac{x^2+(a+b+c)x}{2abc}.
\end{align*}
It is trivial to check that for all integers $k\ge 1$ and $j\ge 0$,
\begin{align*}
f(k+2j)<g(k+2j+1)\le f(k+2j+1)<g(k+2j+2)\le f(k+2j+2).
\end{align*}
For any integer $n\ge f(k)$, there exists a unique integer $l\ge 0$ such that
\begin{align*}
f(k+2l)\le n< f(k+2l+2).
\end{align*}

By Lemma \ref{b-lemma}, we rewrite \eqref{a-gamma} as
\begin{align}
\sum_{n= 0}^{\infty}\gamma_{a,b,c,A,B}^k (n)q^n=\sum_{n=0}^{\infty}\left(F(n)+\beta(n)\right)q^n\sum_{j\not\in [-k,k-1]} (-1)^{j+k}q^{Aj^2+Bj}.\label{c-1}
\end{align}

Next, we shall distinguish nine cases to prove Theorem \ref{t-1}.

{\noindent\bf Case 1} $n<f(k)$. By \eqref{c-1}, it is easy to see that $\gamma_{a,b,c,A,B}^k (n)=0$.
\\[7pt]
{\noindent\bf Case 2} $f(k)\le n<g(k+1)$. By \eqref{c-1}, we have $\gamma_{a,b,c,A,B}^k (n)=\alpha\left(n-f(k)\right)\ge 0$.
\\[7pt]
{\noindent\bf Case 3} $g(k+1)\le n<f(k+1)$. We have
\begin{align*}
\gamma_{a,b,c,A,B}^k (n)&\ge F\left(n-f(k)\right)-F\left(n-g(k+1)\right)-2D_{a,b,c}\\[7pt]
&=\frac{(2k+1)(A-B)n}{abc}+C,
\end{align*}
where $(2k+1)(A-B)/(abc)>0$ and $C$ is independent of $n$.
It follows that
\begin{align*}
\gamma_{a,b,c,A,B}^k (n)
&\ge F\left(g(k+1)-f(k)\right)-2D_{a,b,c}\\[7pt]
&=\frac{2(A-B)^2k^2}{abc}+\frac{(A-B)(2A-2B+a+b+c)k}{abc}\\[7pt]
&+\frac{(A-B)(A-B+a+b+c)}{2abc}-2D_{a,b,c}\\[7pt]
&=H_1,
\end{align*}
where $2(A-B)^2/(abc)>0$.
Then $\gamma_{a,b,c,A,B}^k (n)\ge H_1\ge 0$ for $g(k+1)\le n<f(k+1)$ with $k\ge \mathcal{T}_k\{H_1\}$.
\\[7pt]
{\noindent\bf Case 4} $f(k+1)\le n<g(k+2)$. We have
\begin{align*}
\gamma_{a,b,c,A,B}^k (n)&\ge F\left(n-f(k)\right)-F\left(n-g(k+1)\right)-F\left(n-f(k+1)\right)
-3D_{a,b,c}\\[7pt]
&=-\frac{n^2}{2abc}+Qn+C,
\end{align*}
where $-1/(2abc)<0$ and $Q,C$ are independent of $n$.
It follows that
\begin{align*}
\gamma_{a,b,c,A,B}^k (n)\ge \min\left\{-\frac{f(k+1)^2}{2abc}+Qf(k+1)+C,-\frac{g(k+2)^2}{2abc}+Qg(k+2)+C\right\}.
\end{align*}

On one hand, we have
\begin{align*}
&-\frac{f(k+1)^2}{2abc}+Qf(k+1)+C\\[7pt]
&=F\left(f(k+1)-f(k)\right)-F\left(f(k+1)-g(k+1)\right)-3D_{a,b,c}\\[7pt]
&=\frac{2(A-B)(A+B)k^2}{abc}+\frac{(A-B)(2A+4B+a+b+c)k}{abc}\\[7pt]
&+\frac{(A-B)(A+3B+a+b+c)}{2abc}-3D_{a,b,c}\\[7pt]
&=H_2,
\end{align*}
where $2(A-B)(A+B)/(abc)>0$.

On the other hand, we have
\begin{align*}
&-\frac{g(k+2)^2}{2abc}+Qg(k+2)+C\\[7pt]
&=F\left(g(k+2)-f(k)\right)-F\left(g(k+2)-g(k+1)\right)-F\left(g(k+2)-f(k+1)\right)
-3D_{a,b,c}\\[7pt]
&=\frac{4A(A-B)k^2}{abc}+\frac{2(2A+B)(A-B)k}{abc}-\frac{(A-B)(A-3B+a+b+c)}{abc}-3D_{a,b,c}\\[7pt]
&=H_3,
\end{align*}
where $4A(A-B)/(abc)>0$.
Then $\gamma_{a,b,c,A,B}^k (n)\ge \min\{H_2,H_3\}\ge 0$ for $f(k+1)\le n<g(k+2)$ with $k\ge \mathcal{T}_k\left\{H_2,H_3\right\}$.
\\[7pt]
{\noindent\bf Case 5} $g(k+2)\le n<f(k+2)$. We have
\begin{align*}
&\gamma_{a,b,c,A,B}^k (n)\\[7pt]
&\ge F\left(n-f(k)\right)-F\left(n-g(k+1)\right)-F\left(n-f(k+1)\right)+F\left(n-g(k+2)\right)
-4D_{a,b,c}\\[7pt]
&=-\frac{2(A-B)n}{abc}+C,
\end{align*}
where $-2(A-B)/(abc)<0$ and $C$ is independent of $n$. It follows that
\begin{align*}
&\gamma_{a,b,c,A,B}^k (n)\\[7pt]
&\ge F\left(f(k+2)-f(k)\right)-F\left(f(k+2)-g(k+1)\right)\\[7pt]
&-F\left(f(k+2)-f(k+1)\right)+F\left(f(k+2)-g(k+2)\right)-4D_{a,b,c}\\[7pt]
&=\frac{4A(A-B)k^2}{abc}+\frac{2(A-B)(2A-B)k}{abc}-\frac{(A-B)(A+5B+a+b+c)}{abc}-4D_{a,b,c}\\[7pt]
&=H_4,
\end{align*}
where $4A(A-B)/(abc)>0$.
Then $\gamma_{a,b,c,A,B}^k (n)\ge H_4\ge 0$ for $g(k+2)\le n<f(k+2)$ with $k\ge \mathcal{T}_k\{H_4\}$.
\\[7pt]
{\noindent\bf Case 6} $f(k+2l)\le n<g(k+2l+1)$ with $l\ge 1$.
Note that
\begin{align*}
&F\left(n-f(k+2j)\right)-F\left(n-g(k+2j+1)\right)\\[7pt]
&-F\left(n-f(k+2j+1)\right)+F\left(n-g(k+2j+2)\right)\\[7pt]
&=\frac{(A-B)(6Ak^2+24Akj+24Aj^2+12Ak+24Aj+7A-B-a-b-c-2n)}{abc}.
\end{align*}
It follows that
\begin{align}
&\sum_{j=0}^{l-1}\left(F\left(n-f(k+2j)\right)-F\left(n-g(k+2j+1)\right)\right)\notag\\[7pt]
&+\sum_{j=0}^{l-1}\left(-F\left(n-f(k+2j+1)\right)+F\left(n-g(k+2j+2)\right)\right)\notag\\[7pt]
&=\frac{l(A-B)(6Ak^2+12Akl+8Al^2-A-B-a-b-c-2n)}{abc}.
\label{c-2}
\end{align}

By \eqref{c-2}, we have
\begin{align*}
\gamma_{a,b,c,A,B}^k (n)
&\ge \frac{l(A-B)(6Ak^2+12Akl+8Al^2-A-B-a-b-c-2n)}{abc}\\[7pt]
&+F\left(n-f(k+2l)\right)-(4l+1)D_{a,b,c}\\[7pt]
&=Pn^2+Qn+C,
\end{align*}
where $C$ is independent of $n$ and
\begin{align*}
&P=\frac{1}{2abc}>0,\\
&Q=-\frac{2Ak^2+8Akl+8Al^2+4Al+2Bk-a-b-c}{2abc}.
\end{align*}
It follows that
\begin{align*}
\gamma_{a,b,c,A,B}^k (n)
&\ge Pn^2+Qn+C\\[7pt]
&\ge -\frac{Q^2}{4P}+C\\[7pt]
&=\frac{2(A-B)(2Ak-A-B)l^2}{abc}+\frac{(2k+1)(A-B)(2Ak-A-B)l}{abc}\\[7pt]
&-\frac{(a+b+c)^2}{8abc}-(4l+1)D_{a,b,c}\\[7pt]
&=G_1,
\end{align*}
where $2(A-B)(2Ak-A-B)/(abc)>0$.
Then $\gamma_{a,b,c,A,B}^k (n)\ge G_1\ge 0$ for $f(k+2l)\le n<g(k+2l+1)$ with $l\ge \mathcal{T}_l\{G_1\}$.

By Lemma \ref{c-lemma} and the fact $-(a+b+c)^2/(8abc)-D_{a,b,c}<0$, we have $\mathcal{T}_l\{G_1\}=1$ for
\begin{align*}
H_5&=\frac{4A(A-B)k^2}{abc}+\frac{2(A-B)(2A-B)k}{abc}-\frac{24(A-B)(A+B)+(a+b+c)^2}{8abc}-5D_{a,b,c}\\[7pt]
&\ge 0,
\end{align*}
where $4A(A-B)/(abc)>0$.

For $k\ge \mathcal{T}_k\{H_5\}$, we have $H_5\ge 0$, and so $\mathcal{T}_l\{G_1\}=1$.
It follows that $\gamma_{a,b,c,A,B}^k (n)\ge G_1\ge 0$ for $f(k+2l)\le n<g(k+2l+1)$ with $l\ge 1$ and $k\ge \mathcal{T}_k\{H_5\}$.
\\[7pt]
{\noindent\bf Case 7} $g(k+2l+1)\le n<f(k+2l+1)$ with $l\ge 1$. By \eqref{c-2}, we have
\begin{align*}
\gamma_{a,b,c,A,B}^k (n)
&\ge \frac{l(A-B)(6Ak^2+12Akl+8Al^2-A-B-a-b-c-2n)}{abc}\\[7pt]
&+F\left(n-f(k+2l)\right)-F\left(n-g(k+2l+1)\right)-(4l+2)D_{a,b,c}\\[7pt]
&=Qn+C,
\end{align*}
where $C$ is independent of $n$ and
\begin{align*}
Q=\frac{(2k+2l+1)(A-B)}{abc}>0.
\end{align*}
It follows that
\begin{align*}
&\gamma_{a,b,c,A,B}^k (n)\\[7pt]
&\ge Qg(k+2l+1)+C\\[7pt]
&=\frac{4(A-B)(Ak-B)l^2}{abc}+\frac{(A-B)(4Ak^2+4Ak-6Bk+A-3B+a+b+c)l}{abc}\\[7pt]
&+\frac{(2k+1)(A-B)(2Ak-2Bk+A-B+a+b+c)}{2abc}-(4l+2)D_{a,b,c}\\[7pt]
&=G_2,
\end{align*}
where $4(A-B)(Ak-B)/(abc)>0$.
Then $\gamma_{a,b,c,A,B}^k (n)\ge G_2\ge 0$ for $g(k+2l+1)\le n<f(k+2l+1)$ with $l\ge \mathcal{T}_l \{G_2\}$.

By Lemma \ref{c-lemma}, we have $\mathcal{T}_l\{G_2\}=1$ for
\begin{align*}
H_6&=\frac{4A(A-B)k^2}{abc}+\frac{6(A-B)(2A-B)k}{abc}+\frac{(A-B)(A-11B+a+b+c)}{abc}-4D_{a,b,c}\\[7pt]
&\ge 0,
\end{align*}
and
\begin{align*}
H_7&=\frac{2(A-B)(3A-B)k^2}{abc}+\frac{(A-B)(10A-8B+a+b+c)k}{abc}\\[7pt]
&+\frac{3(A-B)(A-5B+a+b+c)}{2abc}-6D_{a,b,c}\\[7pt]
&\ge 0,
\end{align*}
where $4A(A-B)/(abc)>0$ and $2(A-B)(3A-B)/(abc)>0$.

For $k\ge \mathcal{T}_k \{H_6,H_7\}$, we have $H_6,H_7\ge 0$, and so $\mathcal{T}_l\{G_2\}=1$.
It follows that $\gamma_{a,b,c,A,B}^k (n)\ge G_2\ge 0$ for $g(k+2l+1)\le n<f(k+2l+1)$ with $l\ge 1$ and $k\ge \mathcal{T}_k \{H_6,H_7\}$.
\\[7pt]
{\noindent\bf Case 8} $f(k+2l+1)\le n<g(k+2l+2)$ with $l\ge 1$.
By \eqref{c-2}, we have
\begin{align*}
&\gamma_{a,b,c,A,B}^k (n)\\
&\ge \frac{l(A-B)(6Ak^2+12Akl+8Al^2-A-B-a-b-c-2n)}{abc}\\
&+F\left(n-f(k+2l)\right)-F\left(n-g(k+2l+1)\right)-F\left(n-f(k+2l+1)\right)-(4l+3)D_{a,b,c}\\
&=Pn^2+Qn+C,
\end{align*}
where $Q$ and $C$ are independent of $n$ and
\begin{align*}
P=-\frac{1}{2abc}<0.
\end{align*}
It follows that
\begin{align*}
&\gamma_{a,b,c,A,B}^k (n)\\[7pt]
&\ge \min\left\{Pf(k+2l+1)^2+Qf(k+2l+1)+C,Pg(k+2l+2)^2+Qg(k+2l+2)+C\right\}.
\end{align*}

On one hand, we have
\begin{align*}
&Pf(k+2l+1)^2+Qf(k+2l+1)+C\\[7pt]
&=\frac{4(A-B)(Ak+B)l^2}{abc}+\frac{(A-B)(4Ak^2+4Ak+6Bk+A+5B+a+b+c)l}{abc}\\[7pt]
&+\frac{(2k+1)(A-B)(2Ak+2Bk+A+3B+a+b+c)}{2abc}-(4l+3)D_{a,b,c}\\[7pt]
&=G_3,
\end{align*}
where $4(A-B)(Ak+B)/(abc)>0$.

On the other hand, we have
\begin{align*}
&Pg(k+2l+2)^2+Qg(k+2l+2)+C\\[7pt]
&=\frac{4(A-B)(Ak+B)l^2}{abc}+\frac{(A-B)(4Ak^2+8Ak+2Bk-A+7B-a-b-c)l}{abc}\\[7pt]
&+\frac{(A-B)(4Ak^2+4Ak+2Bk-A+3B-a-b-c)}{abc}-(4l+3)D_{a,b,c}\\[7pt]
&=G_4,
\end{align*}
where $4(A-B)(Ak+B)/(abc)>0$.

It follows that $\gamma_{a,b,c,A,B}^k (n)\ge \min\{G_3,G_4\}\ge 0$ for $f(k+2l+1)\le n<g(k+2l+2)$ with $l\ge \mathcal{T}_l\{G_3,G_4\}$.

By Lemma \ref{c-lemma}, we have $\mathcal{T}_l\{G_3\}=1$ for
\begin{align*}
H_8&=\frac{4A(A-B)k^2}{abc}+\frac{6(A-B)(2A+B)k}{abc}+\frac{(A-B)(A+13B+a+b+c)}{abc}-4D_{a,b,c}\\[7pt]
&\ge 0,\\[7pt]
H_9&=\frac{2(A-B)(3A+B)k^2}{abc}+\frac{(A-B)(10A+a+b+c+10B)k}{abc}\\[7pt]
&+\frac{3(A-B)(A+7B+a+b+c)}{2abc}-7D_{a,b,c}\\[7pt]
&\ge 0,
\end{align*}
and $\mathcal{T}_l\{G_4\}=1$ for
\begin{align*}
H_{10}&=\frac{4A(A-B)k^2}{abc}+\frac{2(A-B)(8A+B)k}{abc}+\frac{(A-B)(A-15B+a+b+c)}{abc}
-4D_{a,b,c}\\[7pt]
&\ge 0,\\[7pt]
H_{11}&=\frac{8A(A-B)k^2}{abc}+\frac{4(A-B)(4A+B)k}{abc}-\frac{2(A-B)(A-7B+a+b+c)}{abc}-7D_{a,b,c}\\[7pt]
&\ge 0,
\end{align*}
where $4A(A-B)/(abc)>0,2(A-B)(3A+B)/(abc)>0$ and $8A(A-B)/(abc)>0$.

For $k\ge \mathcal{T}_k\{H_8,H_9,H_{10},H_{11}\}$, we have $H_8,H_9,H_{10},H_{11}\ge 0$, and so
$\mathcal{T}_l\{G_3\}=\mathcal{T}_l\{G_4\}=1$.
It follows that $\gamma_{a,b,c,A,B}^k (n)\ge \min\{G_3,G_4\}\ge 0$ for $f(k+2l+1)\le n<g(k+2l+2)$ with $l\ge 1$ and $k\ge \mathcal{T}_k\{H_8,H_9,H_{10},H_{11}\}$.
\\[7pt]
{\noindent\bf Case 9} $g(k+2l+2)\le n<f(k+2l+2)$ with $l\ge 1$.
By \eqref{c-2}, we have
\begin{align*}
\gamma_{a,b,c,A,B}^k (n)
&\ge \frac{l(A-B)(6Ak^2+12Akl+8Al^2-A-B-a-b-c-2n)}{abc}\\[7pt]
&+F\left(n-f(k+2l)\right)-F\left(n-g(k+2l+1)\right)-F\left(n-f(k+2l+1)\right)\\[7pt]
&+F\left(n-g(k+2l+2)\right)-(4l+4)D_{a,b,c}\\[7pt]
&=Qn+C,
\end{align*}
where $C$ is independent of $n$ and
\begin{align*}
Q=-\frac{2(l+1)(A-B)}{abc}<0.
\end{align*}
It follows that
\begin{align*}
&\gamma_{a,b,c,A,B}^k (n)\\[7pt]
&\ge Qf(k+2l+2)+C\\[7pt]
&=\frac{4(A-B)(Ak-B)l^2}{abc}+\frac{(A-B)(4Ak^2+8Ak-2Bk-A-9B-a-b-c)l}{abc}\\[7pt]
&+\frac{(A-B)(4Ak^2+4Ak-2Bk-A-5B-a-b-c)}{abc}-(4l+4)D_{a,b,c}\\[7pt]
&=G_5,
\end{align*}
where $4(A-B)(Ak-B)/(abc)>0$.
Then $\gamma_{a,b,c,A,B}^k (n)\ge G_5\ge 0$ for $g(k+2l+2)\le n<f(k+2l+2)$ with $l\ge \mathcal{T}_l\{G_5\}$.

By Lemma \ref{c-lemma}, we have $\mathcal{T}_l\{G_5\}=1$ for
\begin{align*}
H_{12}&=\frac{4A(A-B)k^2}{abc}+\frac{2(A-B)(8A-B)k}{abc}-\frac{(A-B)(A+17B+a+b+c)}{abc}-4D_{a,b,c}\\[7pt]
&\ge 0,\\[7pt]
H_{13}&=\frac{8A(A-B)k^2}{abc}+\frac{4(A-B)(4A-B)k}{abc}-\frac{2(A-B)(A+9B+a+b+c)}{abc}-8D_{a,b,c}\\[7pt]
&\ge 0,
\end{align*}
where $4A(A-B)/(abc)>0$ and $8A(A-B)/(abc)>0$.

For $k\ge \mathcal{T}_k\{H_{12},H_{13}\}$, we have $H_{12},H_{13}\ge 0$, and so $\mathcal{T}_l\{G_5\}=1$.
It follows that $\gamma_{a,b,c,A,B}^k (n)\ge G_5\ge 0$ for $g(k+2l+2)\le n<f(k+2l+2)$ with $l\ge 1$ and $k\ge \mathcal{T}_k\{H_{12},H_{13}\}$.

Combining the above nine cases, we conclude that $\gamma_{a,b,c,A,B}^k (n)\ge 0$ for all $n\ge 0$ with $k\ge \lceil\mathcal{T}_{k}\left\{H_1,H_2,\cdots,H_{13}\right\}\rceil$, and $\gamma_{a,b,c,A,B}^k (n)\ge 0$ for all $n\ge f(k+2l)$ with $l= \lceil\mathcal{T}_{l}\left\{G_1,G_2,\cdots,G_5\right\}\rceil$. This completes the proof of Theorem \ref{t-1}.

\section{Proof of \eqref{merca-conj-1} and \eqref{merca-conj-2}}
In order to prove \eqref{merca-conj-1} and \eqref{merca-conj-2}, we require two preliminary results.
\begin{lem}
We have
\begin{align}
\frac{\prod_{n=1}^{\infty}(1-q^{N_{2n}})}{(q;q)_{\infty}}=\frac{P(q)}{(1-q)(1-q^4)(1-q^5)},\label{e-1}
\end{align}
where $P(q)\in \mathbb{N}[[q]]$.
\end{lem}
{\noindent\it Proof.}
Note that
\begin{align*}
\frac{\prod_{n=1}^{\infty}(1-q^{N_{2n}})}{(q;q)_{\infty}}
&=\frac{\prod_{n=1}^{\infty}(1-q^{N_{4n-2}})}{(q;q^2)_{\infty}}
\cdot \frac{\prod_{n=1}^{\infty}(1-q^{N_{16n-8}})}{(q^4;q^8)_{\infty}}\\[5pt]
&\times \frac{\prod_{n=1}^{\infty}(1-q^{N_{16n}})}{(q^8;q^8)_{\infty}}
\cdot \frac{\prod_{n=1}^{\infty}(1-q^{N_{16n-12}})}{(q^2;q^8)_{\infty}}
\cdot \frac{\prod_{n=1}^{\infty}(1-q^{N_{16n-4}})}{(q^6;q^8)_{\infty}},
\end{align*}
and
\begin{align*}
&N_{4n-2}=3(2n-1),\\[5pt]
&N_{16n-8}=5(8n-4),\\[5pt]
&N_{16n}=8n(6+\nu_2(n)),\\[5pt]
&N_{16n-12}=4(8n-6),\\[5pt]
&N_{16n-4}=4(8n-2).
\end{align*}
Thus, we have
\begin{align}
&\frac{\prod_{n=1}^{\infty}(1-q^{N_{2n}})}{(q;q)_{\infty}}\notag\\[5pt]
&=\frac{\prod_{n=1}^{\infty}(1-q^{3(2n-1)})}{\prod_{n=1}^{\infty}(1-q^{2n-1})}
\cdot \frac{\prod_{n=1}^{\infty}(1-q^{5(8n-4)})}{\prod_{n=1}^{\infty}(1-q^{8n-4})}\notag\\[5pt]
&\times \frac{\prod_{n=1}^{\infty}(1-q^{8n(6+\nu_2(n))})}{\prod_{n=1}^{\infty}(1-q^{8n})}
\cdot \frac{\prod_{n=1}^{\infty}(1-q^{4(8n-6)})}{\prod_{n=1}^{\infty}(1-q^{8n-6})}
\cdot \frac{\prod_{n=1}^{\infty}(1-q^{4(8n-2)})}{\prod_{n=1}^{\infty}(1-q^{8n-2})}.\label{e-2}
\end{align}

For all positive integers $n$, we have
\begin{align*}
&3(2n-1)\equiv 1\pmod{2},\\
&5(8n-4)\equiv 4\pmod{8}.
\end{align*}
It follows that
\begin{align}
\frac{\prod_{n=1}^{\infty}(1-q^{3(2n-1)})}{\prod_{n=1}^{\infty}(1-q^{2n-1})}
\cdot \frac{\prod_{n=1}^{\infty}(1-q^{5(8n-4)})}{\prod_{n=1}^{\infty}(1-q^{8n-4})}
=\frac{P_1(q)}{(1-q)(1-q^4)(1-q^5)},\label{e-3}
\end{align}
where
\begin{align*}
P_1(q)=\frac{\prod_{n=2}^{\infty}(1-q^{3(2n-1)})}{\prod_{n=4}^{\infty}(1-q^{2n-1})}
\cdot \frac{\prod_{n=1}^{\infty}(1-q^{5(8n-4)})}{\prod_{n=2}^{\infty}(1-q^{8n-4})}\in \mathbb{N}[[q]].
\end{align*}

For all positive integers $n$, we have
\begin{align*}
&\frac{1-q^{8n(6+\nu_2(n))})}{1-q^{8n}}\in \mathbb{N}[[q]],\\[5pt]
&\frac{1-q^{4(8n-6)}}{1-q^{8n-6}}\in \mathbb{N}[[q]],\\[5pt]
&\frac{1-q^{4(8n-2)}}{1-q^{8n-2}}\in \mathbb{N}[[q]],
\end{align*}
and so
\begin{align}
P_2(q)=\frac{\prod_{n=1}^{\infty}(1-q^{8n(6+\nu_2(n))})}{\prod_{n=1}^{\infty}(1-q^{8n})}
\cdot \frac{\prod_{n=1}^{\infty}(1-q^{4(8n-6)})}{\prod_{n=1}^{\infty}(1-q^{8n-6})}
\cdot \frac{\prod_{n=1}^{\infty}(1-q^{4(8n-2)})}{\prod_{n=1}^{\infty}(1-q^{8n-2})}\in \mathbb{N}[[q]].
\label{e-4}
\end{align}

Finally, combining \eqref{e-2}--\eqref{e-4} gives
\begin{align*}
\frac{\prod_{n=1}^{\infty}(1-q^{N_{2n}})}{(q;q)_{\infty}}=\frac{P_1(q)P_2(q)}{(1-q)(1-q^4)(1-q^5)},
\end{align*}
where $P_1(q)P_2(q)\in \mathbb{N}[[q]]$.
\qed

\begin{lem}
For all positive integers $k$, we have
\begin{align}
\frac{1}{1-q}\sum_{j\not\in[-k,k]}^{\infty}(-1)^{j+k-1}q^{j(3j+1)/2} \in \mathbb{N}[[q]].\label{e-new}
\end{align}
\end{lem}
{\noindent\it Proof.}
Let $f(j)=j(3j+1)/2$ and $g(j)=j(3j-1)/2$. It is clear that for any positive integer $j$,
\begin{align*}
g(j)<f(j)<g(j+1)<f(j+1).
\end{align*}
For any integer $n\ge g(k+1)$, there exists a unique nonnegative integer $l$ such that
\begin{align*}
g(k+2l+1)\le n<g(k+2l+3).
\end{align*}

Let
\begin{align*}
\frac{1}{1-q}\sum_{j\not\in[-k,k]}^{\infty}(-1)^{j+k-1}q^{j(3j+1)/2}
=\sum_{n=0}^{\infty}\gamma'(n)q^n.
\end{align*}
Note that $1/(1-q)=1+q+q^2+q^3+\cdots$. For $n<g(k+1)$, we have $\gamma'(n)=0$.
For $g(k+2l+1)\le n<f(k+2l+1)$ with $l\ge 0$, we have $\gamma'(n)=1$. For $f(k+2l+1)\le n<g(k+2l+2)$ with $l\ge 0$, we have $\gamma'(n)=2$. For $g(k+2l+2)\le n<f(k+2l+2)$ with $l\ge 0$, we have $\gamma'(n)=1$.
For $f(k+2l+2)\le n<g(k+2l+3)$ with $l\ge 0$, we have $\gamma'(n)=0$. Then $\gamma'(n)\ge 0$ for all integers $n\ge 0$.
\qed

Now we are ready to prove \eqref{merca-conj-1} and \eqref{merca-conj-2}.

{\noindent\it Proof of \eqref{merca-conj-1}}.
By \eqref{a-euler}, we have
\begin{align}
&(-1)^k\left(1-\frac{1}{(q;q)_{\infty}}\sum_{j=1-k}^k (-1)^jq^{j(3j-1)/2}\right)\prod_{n=1}^{\infty}(1-q^{N_{2n}})\notag\\[5pt]
&= (-1)^k\left(1-\frac{1}{(q;q)_{\infty}}\left((q;q)_{\infty}-\sum_{j\not\in [1-k,k]}(-1)^jq^{j(3j-1)/2}\right)\right)\prod_{n=1}^{\infty}(1-q^{N_{2n}})\notag\\[5pt]
&=\frac{\prod_{n=1}^{\infty}(1-q^{N_{2n}})}{(q;q)_{\infty}}\sum_{j\not\in [1-k,k]}(-1)^{j+k}q^{j(3j-1)/2}\notag\\[5pt]
&=\frac{\prod_{n=1}^{\infty}(1-q^{N_{2n}})}{(q;q)_{\infty}}\sum_{j\not\in [-k,k-1]}(-1)^{j+k}q^{j(3j+1)/2}.\label{e-5}
\end{align}
It follows from \eqref{e-1} and \eqref{e-5} that
\begin{align}
&(-1)^k\left(1-\frac{1}{(q;q)_{\infty}}\sum_{j=1-k}^k (-1)^jq^{j(3j-1)/2}\right)\prod_{n=1}^{\infty}(1-q^{N_{2n}})\notag\\[5pt]
&=\frac{P(q)}{(1-q)(1-q^4)(1-q^5)}\sum_{j\not\in [-k,k-1]}(-1)^{j+k}q^{j(3j+1)/2},\label{e-6}
\end{align}
where $P(q)\in \mathbb{N}[[q]]$.

By Corollary \ref{cor-3}, we have
\begin{align}
\frac{1}{(1-q)(1-q^4)(1-q^5)}\sum_{j\not\in [-k,k-1]}(-1)^{j+k}q^{j(3j+1)/2}\in \mathbb{N}[[q]].
\label{e-7}
\end{align}
Combining \eqref{e-6} and \eqref{e-7}, we complete the proof of \eqref{merca-conj-1}.
\qed

{\noindent\it Proof of \eqref{merca-conj-2}.}
By \eqref{a-euler}, we have
\begin{align}
&(-1)^{k-1}\left(1-\frac{1}{(q;q)_{\infty}}\sum_{j=-k}^k (-1)^jq^{j(3j-1)/2}\right)\prod_{n=1}^{\infty}(1-q^{N_{2n}})\notag\\[5pt]
&= (-1)^{k-1}\left(1-\frac{1}{(q;q)_{\infty}}\left((q;q)_{\infty}-\sum_{j\not\in [-k,k]}(-1)^jq^{j(3j-1)/2}\right)\right)\prod_{n=1}^{\infty}(1-q^{N_{2n}})\notag\\[5pt]
&=\frac{\prod_{n=1}^{\infty}(1-q^{N_{2n}})}{(q;q)_{\infty}}\sum_{j\not\in [-k,k]}(-1)^{j+k-1}q^{j(3j-1)/2}\notag\\[5pt]
&=\frac{\prod_{n=1}^{\infty}(1-q^{N_{2n}})}{(q;q)_{\infty}}\sum_{j\not\in [-k,k]}(-1)^{j+k-1}q^{j(3j+1)/2}.\label{e-8}
\end{align}
It follows from \eqref{e-1} and \eqref{e-8} that
\begin{align}
&(-1)^{k-1}\left(1-\frac{1}{(q;q)_{\infty}}\sum_{j=-k}^k (-1)^jq^{j(3j-1)/2}\right)\prod_{n=1}^{\infty}(1-q^{N_{2n}})\notag\\[5pt]
&=\frac{P(q)}{(1-q)(1-q^4)(1-q^5)}\sum_{j\not\in [-k,k]}(-1)^{j+k-1}q^{j(3j+1)/2},\label{e-9}
\end{align}
where $P(q)\in \mathbb{N}[[q]]$.
Combining \eqref{e-new} and \eqref{e-9}, we complete the proof of \eqref{merca-conj-2}.
\qed

\vskip 5mm \noindent{\bf Acknowledgments.} This work was supported by the National Natural Science Foundation of China (grant 12171370).

\end{document}